\documentclass[10pt,a4paper,draft]{article}
\usepackage{bbm}
\usepackage{pifont}
\usepackage{bbding}
\pdfoutput=1
\usepackage{amsmath, esint}
\usepackage{mathrsfs}
\usepackage{amsfonts}
\usepackage{amssymb}
\usepackage{amsfonts,amssymb,amsmath,indentfirst,cases,amsthm}
\usepackage{epsfig}
\usepackage[numbers,sort&compress]{natbib}
\def\ess~inf{\mathop{\rm ess~inf}}

\numberwithin{equation}{section}

\newenvironment{key words}{\emph{\texttt{Keywords}}\mbox{  }}{ }
\newtheorem{theorem}{Theorem}[section]

\newtheorem{remark}[theorem]{Remark}

\theoremstyle{remark}

\theoremstyle{plain}

\makeatletter

\newcommand{\Rmnum}[1]{\expandafter\@slowromancap\romannumeral #1@}
\makeatother
\begin{document}

\title{\textbf{A Proof of Riemann Hypothesis}}
\author{Pengcheng Niu$^1$, Junli Zhang$^2$\thanks{Corresponding
author's E-mail:jlzhang2020@163.com}, Zhijun Qiao$^3$\\
\small{1. School of Mathematics and Statistics, Northwestern Polytechnical
University,}\\ \small{ Xi'an, Shaanxi, 710129, P. R. China}\\
\small{2. School of Mathematics and Data Science, Shaanxi University of Science and Technology,}\\
\small{ Xi'an, Shaanxi, 710021, P. R. China}\\
\small{3. School of Mathematical and Statistical Sciences, University of Texas Rio Grande Valley,}\\
\small{ Edinburg, TX 78539, USA}}
\date{} \maketitle

% ----------------------------------------------------------------
\maketitle {\bf Abstract.}\
Let $\Xi(t)$  be a function relating to the Riemann zeta function $\zeta (s)$  with $s = \frac{1}
{2} + it$. In this paper, we construct a function $v$  containing $t$  and $\Xi(t)$, and prove that $v$  satisfies a nonadjoint boundary value problem to a nonsingular differential equation if $t$ is any nontrivial zero of $\Xi(t)$. Inspecting properties of $v$ and using known results of nontrivial zeros of $\zeta (s)$, we derive that nontrivial zeros of $\zeta (s)$  all have real part equal to $\frac{1}
{2}$, which concludes that Riemann Hypothesis is true.\\

\textbf{keywords.} Riemann hypothesis, Riemann zeta function


\def\Xint#1{\mathchoice
    {\XXint\displaystyle\textstyle{#1}}%
    {\XXint\textstyle\scriptstyle{#1}}%
    {\XXint\scriptstyle\scriptscriptstyle{#1}}%
    {\XXint\scriptscriptstyle\scriptscriptstyle{#1}}%
    \!\int}
\def\XXint#1#2#3{{\setbox0=\hbox{$#1{#2#3}{\int}$}
    \vcenter{\hbox{$#2#3$}}\kern-.5\wd0}}
\def\dashint{\Xint-}
%---------------------------------------------------------------------------------
\section{\textbf{Introduction}\label{Section 1}}

Riemann showed in \cite{R} that the function
\begin{equation}\label{eq1.1}
\zeta (s) = \sum\limits_{n = 1}^\infty  {\frac{1}
{{{n^s}}}}
\end{equation}
defined on the domain $\operatorname{Re} s > 1$  can be analytically extended to the whole complex plane $\mathbb{C}$  with a unique pole $s=1$. Riemann Hypothesis (RH) means that nontrivial zeros of function $\zeta (s)$ all have real part equal to $\frac{1}
{2}$. One has known that nontrivial zeros of $\zeta (s)$  all locate in the domain $0 < \operatorname{Re} (s) < 1$ in $\mathbb{C}$. Let
\begin{equation}\label{eq1.2}
\xi (s) = \frac{1}
{2}s(s - 1){\pi ^{ - s/2}}\Gamma (\frac{s}
{2})\zeta (s),
\end{equation}
where $\Gamma (s)$ is the Gamma function, then  $\xi (s)$ is an entire function on $\mathbb{C}$ and zeros of $\xi (s)$  coincide with nontrivial ones of $\zeta (s)$ (Apostol \cite{A}, Pan and Pan \cite{PP}, Lu \cite{L}). The function $\xi (s)$ owns infinitely many zeros; these zeros are symmetric with respect to the real axis, the line $\operatorname{Re} s{\text{ = }}\frac{1}
{2}$  and the point $( \frac{1}{2}, 0)$ in $\mathbb{C},$ and real parts of zeros are in $(0, 1)$, see \cite{A} and \cite{PP}.

If $s{\text{ = }}\operatorname{Re} s + i\operatorname{Im} s$  is a zero of $\xi (s)$, then $\left| {\operatorname{Im} s} \right| > 6$ (\cite{PP}). We denote
\begin{equation}\label{eq1.3}
s: = \frac{1}
{2} + it,
\end{equation}
where
\[\operatorname{Re} t = \operatorname{Im} s \in \mathbb{R} ~\text{~and~} \operatorname{Im} t = \frac{1}
{2} - \operatorname{Re} s \in ( - \frac{1}
{2},\frac{1}
{2}),\]
then
\begin{equation}\label{eq1.4}
\left| {\operatorname{Re} t} \right| > 6.
\end{equation}
Letting
\begin{equation}\label{eq1.5}
\Xi(t) = \xi (\frac{1}
{2} + it),
\end{equation}
RH is stated by
\begin{equation}\label{eq1.6}
\Xi(t) = 0 \Rightarrow t \in \mathbb{R}.
\end{equation}

Following \cite{R}, the function $\Xi(t)$ can be expressed as the Fourier cosine integral
\[\Xi(t) = 4\int_1^\infty  {\frac{{d({x^{\frac{3}
{2}}}\psi '(x))}}
{{dx}}} {x^{ - \frac{1}
{4}}}\cos (\frac{t}
{2}\log x)dx,\]
where $\psi (x) = \sum\limits_{n = 1}^\infty  {\exp \{  - {n^2}\pi x\} } , x \in (1,\infty ).$  With the evident change of variables $y = \frac{1}
{2}\log x$, it becomes
\[\Xi(t) = 2\int_0^\infty  {\cos (tx)} \Phi (x)dx,\]
where
\[\Phi (x) = 2\pi \exp \{ \frac{{5x}}
{2}\} \sum\limits_{n = 1}^\infty  {(2\pi \exp \{ 2x\} {n^2} - 3)} {n^2}\exp \{  - {n^2}\pi \exp \{ 2x\} \}. \]

There were numerous mathematicians studying RH and many important results appeared. Here we do not list detailed literature, but point out that P\'{o}lya in \cite{P} deduced that a function relating to a Sturm-Liouville type operator and being similar to $\Xi(t)$ only has real zeros. A relative conjecture to RH is the Hilbert-P\'{o}lya conjecture (\cite{L}) which says that nontrivial zeros of $\zeta (s)$  are corresponding to eigenvalues of some Hermite operator. Concretely, the Hilbert-P\'{o}lya conjecture illustrates that if the nontrivial zero of $\zeta (s)$  is written as the form $s = \frac{1}
{2} + it$, then  $t$ is corresponding to the eigenvalue of some Hermite operator; using the known fact that eigenvalues of a Hermite operator are real, it implies that $t$  is real and so that nontrivial zeros of $\zeta (s)$ all have real part equal to $\frac{1}
{2}$, which proves Riemann Hypothesis.

M. R. Pistorius proposed an idea in an unpublished paper that one may construct a function implicitly involving $\Xi(t)$; derive a boundary value problem for the function at the zero  $t$ of $\Xi(t)$ and use classic Sturm-Liouville to deduce that $t$  is real.

Recently Bender, Brody and M\"{u}ller \cite{BBM} found a nonadjoint Hamiltonian and indicated that one may prove RH by investigating the reality of eigenvalues for the nonadjoint Hamiltonian.

The main result of the paper is

\textbf{Theorem 1} It follows that (\ref{eq1.6}) is true and so RH holds.

To prove (\ref{eq1.6}), one needs to verify that the zero $t = {t_1} + i{t_2}$ of $\Xi(t)$, ${t_1}{\text{ = }}\operatorname{Re} t, {t_2}{\text{ = Im}}t,$  satisfies ${t_2} = 0$. For certain, we construct a function $v$  containing  $t$ and $\Xi(t)$, where $t$  is a complex parameter; letting that $t$  satisfies $\Xi(t)=0$  and supposing ${t_2} \ne 0$, one can inspect properties of $v$ by using known facts of nontrivial zeros of $\zeta (s)$ and derive that $t^2$ is real, which will lead to a contradiction.

One of difficulties is how to construct a suitable function $v$. We adopt hyperbolic functions in the expression of  $v$, which can be used to infer that $v$ satisfies a nonadjoint boundary value problem to a nonsingular differential equation at the zeros of $\Xi(t)$. Also, to assert that $t^2$  is real, we add a quadratic function and a small parameter $\varepsilon$  in $v$. Besides, how to deal with the nonadjoint boundary value problem is a new difficulty. We multiply the conjugate function of $v$ to the equation and continue to work. The last difficulty is to determine $b_1$ such that $Q(b_1)>0$, see (\ref{eq2.20}) below. We are inspired by numeric analysis in the selection of $b_1$.

The proof of Theorem 1 is given in Section 2. Some formulas used in Section 2 are proved in Section 3.

%%---------------------------------------------------------------------------------
\section{Proof of Theorem 1}\label{Section 2}
%%---------------------------------------------------------------------------------
If $t = {t_1} + i{t_2}$ satisfies $\Xi(t)=0$, where ${t_1} \in \mathbb{R}$  and ${t_2} \in ( - \frac{1}
{2},\frac{1}
{2})$, then it is easy to see that  $t$ is symmetric with respect to $t_1$  axis, $t_2$ axis and the origin, and $\left| {{t_1}} \right| > 6$.

\subsection{The construction of function $v$}

Introduce a function
\begin{equation}\label{eq2.1}
v(y;t,\varepsilon ): = ch[t(y - \frac{b}
{2})] + \frac{1}
{{2\varepsilon }}t{(y - \frac{b}
{2})^2} + \Xi(t)y,
\end{equation}
where $y \in [0,b], ~b$  is positive and real, $t \in \mathbb{C}$ is a complex parameter and $\varepsilon$  is a positive parameter. Clearly, $v(y;t,\varepsilon )$  is  ${C^\infty }$ in $y$.

It knows from (\ref{eq2.1}) that
\[v(0;t,\varepsilon ) = ch(t \frac{{ - b}}
{2}) + \frac{1}
{{2\varepsilon }}t {(\frac{{ - b}}
{2})^2},\]
\[v(b;t,\varepsilon ) = ch(t  \frac{b}
{2}) + \frac{1}
{{2\varepsilon }}t{(\frac{b}
{2})^2} + b \Xi(t).\]
Differentiating (\ref{eq2.1}) in $y$, we have
\[v'(y;t,\varepsilon ) = tsh[t(y - \frac{b}
{2})] + \frac{1}
{\varepsilon }t(y - \frac{b}
{2}) + \Xi(t),\]
and
\[v'(0;t,\varepsilon ) = tsh[t(\frac{{ - b}}
{2})] + \frac{1}
{\varepsilon }t(\frac{{ - b}}
{2}) + \Xi(t),\]
\[v'(b;t,\varepsilon ) = tsh(t\frac{b}
{2}) + \frac{1}
{\varepsilon }t(\frac{b}
{2}) + \Xi(t).\]
Obviously,
\[v''(y;t,\varepsilon ) = {t^2}ch[t(y - \frac{b}
{2})] + \frac{1}
{\varepsilon }t = {t^2}v(y;t,\varepsilon ) - \frac{1}
{{2\varepsilon }}{t^3}{(y - \frac{b}
{2})^2} - {t^2} \Xi(t) y + \frac{1}
{\varepsilon }t.\]

Now let  $t = {t_1} + i{t_2}$ satisfy
\begin{equation}\label{eq2.2}
\Xi(t) = 0
\end{equation}
(therefore $t$, $t_1$ and $t_2$ are fixed), then we have from above that

\begin{equation}\label{eq2.3}
v(y;\varepsilon ): = ch[t(y - \frac{b}
{2})] + \frac{1}
{{2\varepsilon }}t{(y - \frac{b}
{2})^2},
\end{equation}
\begin{equation}\label{eq2.4}
v(0;\varepsilon ) = v(b;\varepsilon ) = ch(t  \frac{b}
{2}) + \frac{1}
{{2\varepsilon }}t  {(\frac{b}
{2})^2},
\end{equation}
\[v'(y;\varepsilon ) = tsh[t(y - \frac{b}
{2})] + \frac{1}
{\varepsilon }t(y - \frac{b}
{2}),\]
\begin{equation}\label{eq2.5}
v'(0;\varepsilon ) =  - v'(b;\varepsilon ) =  - tsh(t\frac{b}
{2}) - \frac{1}
{\varepsilon }t\frac{b}
{2},
\end{equation}
\begin{equation}\label{eq2.6}
v''(y;\varepsilon ) = {t^2}v(y;\varepsilon ) - \frac{1}
{{2\varepsilon }}{t^3}{(y - \frac{b}
{2})^2} + \frac{1}
{\varepsilon }t.
\end{equation}

Observe that the equation (\ref{eq2.6}) with the boundary value conditions (\ref{eq2.4}) and (\ref{eq2.5}) form a nonadjoint boundary value problem to a nonsingular differential equation. It can not help us leads to that $t^2$ is real by using the classic Sturm-Liouville theory.

We want to prove $t_2=0$. Let us use the contradiction and assume ${t_2} \ne 0$.
Without loss of generality, one can assume
\[0 < {t_2} < \frac{1}
{2}.\]
Following the symmetry of zeros of $\Xi(t)$ to the  $t_2$ axis, we also let
\[{t_1} > 6.\]

For simplicity, denote $v(y;\varepsilon )$  by $v$, $v'(y;\varepsilon )$ by $v'$,  $v''(y;\varepsilon )$ by $v''$,  $v(0;\varepsilon )$ by $v(0)$, $v(b;\varepsilon )$ by $v(b)$, $v'(0;\varepsilon )$ by  $v'(0)$ and  $v'(b;\varepsilon )$ by $v'(b)$, respectively.

Multiplying $\bar v = ch[\bar t(y - \frac{b}
{2})] + \frac{1}
{{2\varepsilon }}\bar t{(y - \frac{b}
{2})^2}$  (the conjugate function of $v$) on (\ref{eq2.6}) and integrating in $y \in [0,b]$, we have
\[\int_0^b {v''} \bar vdy = \int_0^b {[{t^2}v} \bar v - \frac{1}
{{2\varepsilon }}{t^3}{(y - \frac{b}
{2})^2}\bar v + \frac{1}
{\varepsilon }t\bar v]dy\]
and
\[\int_0^b {v''} \bar vdy = \int_0^b {\bar v} dv' = \bar vv'|_0^b - \int_0^b {v'} \bar v'dy = 2\bar v(b)v'(b) - \int_0^b {v'} \bar v'dy,\]
then
\begin{equation}\label{eq2.7}
{t^2}\int_0^b v \bar vdy - \frac{1}
{{2\varepsilon }}{t^3}\int_0^b {{(y - \frac{b}
{2})^2}} \bar vdy + \frac{1}
{\varepsilon }t\int_0^b {\bar v} dy = 2\bar v(b)v'(b) - \int_0^b {v'} \bar v'dy.
\end{equation}
Let us calculate the second and third terms in the left hand side and the first term in the right hand side of (\ref{eq2.7}), respectively. Using ${y_1} = y - \frac{b}
{2}$ and
\begin{align*}
&~~~~~\int {y_1^2[ch(\bar t{y_1}) + \frac{{\bar t}}
{{2\varepsilon }}y_1^2} ]d{y_1}\\
& = \int {y_1^2\frac{1}
{{\bar t}}dsh(\bar t{y_1}) + \int {\frac{{\bar t}}
{{2\varepsilon }}y_1^4} } d{y_1}\\
& = \frac{1}
{{\bar t}}y_1^2sh(\bar t{y_1}) - \frac{2}
{{{{\bar t}^2}}}{y_1}ch(\bar t{y_1}) + \frac{2}
{{{{\bar t}^3}}}sh(\bar t{y_1}) + \frac{{\bar t}}
{{10\varepsilon }}y_1^5 + constant,
\end{align*}
the second term in the left hand side of (\ref{eq2.7}) is by (\ref{eq2.3}) that

\begin{align}\label{eq2.8}
&~~~~~ - \frac{1}
{{2\varepsilon }}{t^3}\int_0^b {{(y - \frac{b}
{2})^2}\bar v} dy\nonumber\\
& =  - \frac{1}
{{2\varepsilon }}{t^3}\int_0^b {{(y - \frac{b}
{2})^2}\{ ch[\bar t(y - \frac{b}
{2})] + \frac{{\bar t}}
{{2\varepsilon }}{(y - \frac{b}
{2})^2}\} } dy\nonumber\\
& =  - \frac{1}
{{2\varepsilon }}{t^3}\int_{ - \frac{b}
{2}}^{\frac{b}
{2}} {y_1^2[ch(\bar t{y_1}) + \frac{{\bar t}}
{{2\varepsilon }}y_1^2]} dy_1\nonumber\\
& =  - \frac{1}
{\varepsilon }{t^3}[\frac{1}
{{\bar t}}{(\frac{b}
{2})^2}sh(\bar t\frac{b}
{2}) - \frac{2}
{{{{\bar t}^2}}}\frac{b}
{2}ch(\bar t\frac{b}
{2}) + \frac{2}
{{{{\bar t}^3}}}sh(\bar t\frac{b}
{2}) + \frac{{\bar t}}
{{10\varepsilon }}{(\frac{b}
{2})^5}]\nonumber\\
&= - \frac{1}
{\varepsilon }\frac{{{t^4}}}
{{t\bar t}}{(\frac{b}
{2})^2}sh(\bar t\frac{b}
{2}) + \frac{2}
{\varepsilon }\frac{{{t^5}}}
{{{{(t\bar t)}^2}}}\frac{b}
{2}ch(\bar t\frac{b}
{2}) - \frac{2}
{\varepsilon }\frac{{{t^6}}}
{{{{(t\bar t)}^3}}}sh(\bar t\frac{b}
{2}) - \frac{{{t^2}(t\bar t)}}
{{10{\varepsilon ^2}}}{(\frac{b}
{2})^5},
\end{align}
where $t\bar t = t_1^2 + t_2^2$. The third term in the left hand side of (\ref{eq2.7}) satisfies
\begin{align}\label{eq2.9}
\frac{1}
{\varepsilon }t\int_0^b \bar vdy &=  \frac{1}
{\varepsilon }t\int_0^b {\{ ch[\bar t(y - \frac{b}
{2})] + \frac{1}
{{2\varepsilon }}\bar t{(y - \frac{b}
{2})^2}\} dy}\nonumber \\
& = \frac{1}
{\varepsilon }t\{ \frac{1}
{{\bar t}}sh[\bar t(y - \frac{b}
{2})] + \frac{1}
{{2\varepsilon }}\bar t \cdot \frac{1}
{3}{(y - \frac{b}
{2})^3}\} |_0^b\nonumber\\
& = \frac{2}
{\varepsilon }\frac{t}
{{\bar t}}sh(\bar t\frac{b}
{2}) + \frac{1}
{{3{\varepsilon ^2}}}t\bar t{(\frac{b}
{2})^3}\nonumber\\
& = \frac{2}
{\varepsilon }\frac{{{t^2}}}
{{t\bar t}}sh(\bar t\frac{b}
{2}) + \frac{1}
{{3{\varepsilon ^2}}}t\bar t{(\frac{b}
{2})^3}.
\end{align}
The first term in the right hand side of (\ref{eq2.7}) becomes
\begin{align}\label{eq2.10}
2\bar v(b)v'(b) &= 2[ch(\bar t\frac{b}
{2}) + \frac{1}
{{2\varepsilon }}\bar t{(\frac{b}
{2})^2}][tsh(t\frac{b}
{2}) + \frac{1}
{\varepsilon }t\frac{b}
{2}]\nonumber\\
& = 2[tch(\bar t\frac{b}
{2})sh(t\frac{b}
{2}) + \frac{1}
{\varepsilon }t\frac{b}
{2}ch(\bar t\frac{b}
{2}) + \frac{1}
{{2\varepsilon }}t\bar t{(\frac{b}
{2})^2}sh(t\frac{b}
{2}) + \frac{1}
{{2{\varepsilon ^2}}}t\bar t{(\frac{b}
{2})^3}].
\end{align}
Putting (\ref{eq2.8}), (\ref{eq2.9}) and (\ref{eq2.10}) into (\ref{eq2.7}), it yields
\begin{align}\label{eq2.11}
&~~~~~{t^2}\int_0^b {v\bar vdy - \frac{1}
{\varepsilon }\frac{{{t^4}}}
{{t\bar t}}{(\frac{b}
{2})^2}sh(\bar t\frac{b}
{2}) + \frac{2}
{\varepsilon }\frac{{{t^5}}}
{{{{(t\bar t)}^2}}}\frac{b}
{2}ch(\bar t\frac{b}
{2}) - \frac{2}
{\varepsilon }\frac{{{t^6}}}
{{{{(t\bar t)}^3}}}sh(\bar t\frac{b}
{2}) - \frac{{{t^2}t\bar t}}
{{10{\varepsilon ^2}}}{(\frac{b}
{2})^5}} \nonumber\\
&~~~~~ + \frac{2}
{\varepsilon }\frac{{{t^2}}}
{{t\bar t}}sh(\bar t\frac{b}
{2}) + \frac{1}
{{3{\varepsilon ^2}}}t\bar t{(\frac{b}
{2})^3}\nonumber\\
&= - \int_0^b {v'} \bar {v'} dy + 2tch(\bar t\frac{b}
{2})sh(t\frac{b}
{2}) + \frac{2}
{\varepsilon }t\frac{b}
{2}ch(\bar t\frac{b}
{2}) + \frac{1}
{\varepsilon }t\bar t{(\frac{b}
{2})^2}sh(t\frac{b}
{2}) + \frac{1}
{{{\varepsilon ^2}}}t\bar t{(\frac{b}
{2})^3}.
\end{align}

Denoting
\begin{align}\label{eq2.12}
 P(b;\varepsilon )= & - \frac{1}
{\varepsilon }\frac{{{t^4}}}
{{t\bar t}}{(\frac{b}
{2})^2}sh(\bar t\frac{b}
{2}) + \frac{2}
{\varepsilon }\frac{{{t^5}}}
{{{{(t\bar t)}^2}}}\frac{b}
{2}ch(\bar t\frac{b}
{2}) - \frac{2}
{\varepsilon }\frac{{{t^6}}}
{{{{(t\bar t)}^3}}}sh(\bar t\frac{b}
{2}) + \frac{2}
{\varepsilon }\frac{{{t^2}}}
{{t\bar t}}sh(\bar t\frac{b}
{2})\nonumber\\
& - 2tch(\bar t\frac{b}
{2})sh(t\frac{b}
{2}) - \frac{2}
{\varepsilon }t\frac{b}
{2}ch(\bar t\frac{b}
{2}) - \frac{1}
{\varepsilon }t\bar t{(\frac{b}
{2})^2}sh(t\frac{b}
{2}),
\end{align}
then (2.11) can be written as
\begin{equation}\label{eq2.13}
{t^2}\int_0^b {v\bar vdy + P(b;\varepsilon ) - \frac{{{t^2}t\bar t}}
{{10{\varepsilon ^2}}}{(\frac{b}
{2})^5}}  - \frac{2}
{{3{\varepsilon ^2}}}t\bar t{(\frac{b}
{2})^3} + \int_0^b {v'} \bar {v'} dy = 0.
\end{equation}

The next aim is to determine ${b_0} > 0$  and a small ${\varepsilon _0} > 0$  such that
\begin{equation}\label{eq2.14}
\text{the imaginary part of~} P({b_0};{\varepsilon _0})  \text{~vanishes},
\end{equation}
(i.e., $P({b_0};{\varepsilon _0})$ is real) and
\begin{equation}\label{eq2.15}
\int_0^{{b_0}} {v\bar vdy}  - \frac{{t\bar t}}
{{10\varepsilon _0^2}}{(\frac{{{b_0}}}
{2})^5} \ne 0;
\end{equation}
then noting $\int_0^{{b_0}} {v\bar vdy} ( = \int_0^{{b_0}} {{{\left| v \right|}^2}dy} )$, $\frac{{t\bar t}}
{{10\varepsilon _0^2}}{(\frac{{{b_0}}}
{2})^5},$  $ - \frac{2}
{{3\varepsilon _0^2}}t\bar t{(\frac{{{b_0}}}
{2})^3}$  and $\int_0^{{b_0}} {v'} \overline {v'} dy
( = \int_0^{{b_0}} {{{\left| {v'} \right|}^2}} dy )$ in (\ref{eq2.13}) are all real and not equal to zero, we will derive that $t^2$  is real.

\subsection{ The imaginary part $f(b;\varepsilon )$ of $P(b;\varepsilon )$}

To calculate $P(b;\varepsilon )$, recall that for $x$  and $y \in \mathbb{R}$,
\[sh(x \pm iy) = shx\cos y \pm ichx\sin y\]
and
\[ch(x \pm iy) = chx\cos y \pm ishx\sin y,\]
and note
\[{t^2} = ({t_1} + i{t_2})({t_1} + i{t_2}) = t_1^2 - t_2^2 + i \cdot 2{t_1}{t_2},\]
\[{t^4} = {t^2} \cdot {t^2} = {(t_1^2 - t_2^2)^2} - {(2{t_1}{t_2})^2} + i \cdot 2(t_1^2 - t_2^2) \cdot (2{t_1}{t_2}),\]
\begin{align*}
{t^5} &= t \cdot {t^4}\\
& = {t_1}[{(t_1^2 - t_2^2)^2} - {(2{t_1}{t_2})^2}] - {t_2} \cdot 2(t_1^2 - t_2^2) \cdot (2{t_1}{t_2})\\
&~~~ + i\{ {t_1} \cdot 2(t_1^2 - t_2^2) \cdot (2{t_1}{t_2}) + {t_2}[{(t_1^2 - t_2^2)^2} - {(2{t_1}{t_2})^2}]\} ,
\end{align*}
\begin{align*}
{t^6} &= {t^2} \cdot {t^4}\\
& = (t_1^2 - t_2^2)[{(t_1^2 - t_2^2)^2} - {(2{t_1}{t_2})^2}] - (2{t_1}{t_2}) \cdot 2(t_1^2 - t_2^2)(2{t_1}{t_2})\\
&~~~ + i\{ (t_1^2 - t_2^2) \cdot 2(t_1^2 - t_2^2)(2{t_1}{t_2}) + 2{t_1}{t_2}[{(t_1^2 - t_2^2)^2} - {(2{t_1}{t_2})^2}]\} \\
& = (t_1^2 - t_2^2)[{(t_1^2 - t_2^2)^2} - {(2{t_1}{t_2})^2}] - 2(t_1^2 - t_2^2){(2{t_1}{t_2})^2}\\
&~~~ + i \cdot 2{t_1}{t_2}[3{(t_1^2 - t_2^2)^2} - {(2{t_1}{t_2})^2}],
\end{align*}
\[sh(\bar t\frac{b}
{2}) = sh({t_1}\frac{b}
{2} - i{t_2}\frac{b}
{2}) = sh({t_1}\frac{b}
{2})\cos ({t_2}\frac{b}
{2}) - ich({t_1}\frac{b}
{2})\sin ({t_2}\frac{b}
{2}),\]
\[sh(t\frac{b}
{2}) = sh({t_1}\frac{b}
{2} + i{t_2}\frac{b}
{2}) = sh({t_1}\frac{b}
{2})\cos ({t_2}\frac{b}
{2}) + ich({t_1}\frac{b}
{2})\sin ({t_2}\frac{b}
{2}),\]
\[ch(\bar t\frac{b}
{2}) = ch({t_1}\frac{b}
{2} - i{t_2}\frac{b}
{2}) = ch({t_1}\frac{b}
{2})\cos ({t_2}\frac{b}
{2}) - ish({t_1}\frac{b}
{2})\sin ({t_2}\frac{b}
{2}),\]
\begin{align*}
ch(\bar t\frac{b}
{2})sh(t\frac{b}
{2}) &= \frac{1}
{2}({e^{\bar t\frac{b}
{2}}} + {e^{ - \bar t\frac{b}
{2}}}) \cdot \frac{1}
{2}({e^{t\frac{b}
{2}}} - {e^{ - t\frac{b}
{2}}})\\
& = \frac{1}
{4}({e^{(t + \bar t)\frac{b}
{2}}} - {e^{(\bar t - t)\frac{b}
{2}}} + {e^{(t - \bar t)\frac{b}
{2}}} - {e^{( - t - \bar t)\frac{b}
{2}}})\\
& = \frac{1}
{4}({e^{2{t_1}\frac{b}
{2}}} - {e^{ - 2i{t_2}\frac{b}
{2}}} + {e^{2i{t_2}\frac{b}
{2}}} - {e^{ - 2{t_1}\frac{b}
{2}}})\\
& = \frac{1}
{2}[sh({t_1}b) + i\sin ({t_2}b)].
\end{align*}
Inserting these into (\ref{eq2.12}), we have
\begin{align*}
&~~~~~P(b;\varepsilon )\\
& =  - \frac{1}
{\varepsilon }{(\frac{b}
{2})^2}\frac{1}
{{t\bar t}}[{(t_1^2 - t_2^2)^2} - {(2{t_1}{t_2})^2} + i \cdot 2(t_1^2 - t_2^2) \cdot (2{t_1}{t_2})]\\
& ~~~~~~~~~~~~~~~~~\cdot [sh({t_1}\frac{b}
{2})\cos ({t_2}\frac{b}
{2}) - ich({t_1}\frac{b}
{2})\sin ({t_2}\frac{b}
{2})]\\
&~~~ + \frac{2}
{\varepsilon }\frac{b}
{2}\frac{1}
{{{{(t\bar t)}^2}}}\{ {t_1}[{(t_1^2 - t_2^2)^2} - {(2{t_1}{t_2})^2}] - {t_2} \cdot 2(t_1^2 - t_2^2) \cdot (2{t_1}{t_2})\,\\
&~~~~~~~~~~~~~~~~\; + i\{ {t_1} \cdot 2(t_1^2 - t_2^2) \cdot (2{t_1}{t_2}) + {t_2}[{(t_1^2 - t_2^2)^2} - {(2{t_1}{t_2})^2}]\} \} \\
&~~~~~~~~~~~~~~~~~ \cdot [ch({t_1}\frac{b}
{2})\cos ({t_2}\frac{b}
{2}) - ish({t_1}\frac{b}
{2})\sin ({t_2}\frac{b}
{2})]
\end{align*}

\begin{align*}
& ~~~- \frac{2}
{\varepsilon }\frac{1}
{{{{(t\bar t)}^3}}}\{ (t_1^2 - t_2^2)[{(t_1^2 - t_2^2)^2} - {(2{t_1}{t_2})^2}] - 2(t_1^2 - t_2^2){(2{t_1}{t_2})^2}\, + i \cdot 2{t_1}{t_2}[3{(t_1^2 - t_2^2)^2} - {(2{t_1}{t_2})^2}]\} \\
&~~~~~~~~~~~~~~~\cdot [sh({t_1}\frac{b}
{2})\cos ({t_2}\frac{b}
{2}) - ich({t_1}\frac{b}
{2})\sin ({t_2}\frac{b}
{2})]\\
&~~~+ \frac{2}
{\varepsilon }\frac{1}
{{t\bar t}}[t_1^2 - t_2^2 + i \cdot 2{t_1}{t_2}] \cdot [sh({t_1}\frac{b}
{2})\cos ({t_2}\frac{b}
{2}) - ich({t_1}\frac{b}
{2})\sin ({t_2}\frac{b}
{2})]\\
&~~~ - 2({t_1} + i{t_2}) \cdot \frac{1}
{2}[sh({t_1}b) + i\sin ({t_2}b)]\\
&~~~ - \frac{2}
{\varepsilon }\frac{b}
{2}({t_1} + i{t_2}) \cdot [ch({t_1}\frac{b}
{2})\cos ({t_2}\frac{b}
{2}) - ish({t_1}\frac{b}
{2})\sin ({t_2}\frac{b}
{2})]\\
&~~~ - \frac{1}
{\varepsilon }{(\frac{b}
{2})^2}t\bar t[sh({t_1}\frac{b}
{2})\cos ({t_2}\frac{b}
{2}) + ich({t_1}\frac{b}
{2})\sin ({t_2}\frac{b}
{2})].
\end{align*}
Denoting the imaginary part Im$P(b;\varepsilon )$  of  $P(b;\varepsilon )$ by $f(b;\varepsilon )$  ( $f(b;\varepsilon )$ is a real function in $y$ ), it follows
\begin{align}\label{eq2.16}
&~~~~~f(b;\varepsilon )\nonumber\\
& =  - \frac{1}
{\varepsilon }{(\frac{b}
{2})^2}\frac{1}
{{t\bar t}}\{ [{(t_1^2 - t_2^2)^2} - {(2{t_1}{t_2})^2}] \cdot [ - ch({t_1}\frac{b}
{2})\sin ({t_2}\frac{b}
{2})]\nonumber\\
&~~~~~~~~~~~~~~~~~ + 2(t_1^2 - t_2^2) \cdot (2{t_1}{t_2})sh({t_1}\frac{b}
{2})\cos ({t_2}\frac{b}
{2})\}\nonumber\\
&~~~ + \frac{2}
{\varepsilon }\frac{b}
{2}\frac{1}
{{{{(t\bar t)}^2}}}\{ \{ {t_1}[{(t_1^2 - t_2^2)^2} - {(2{t_1}{t_2})^2}] - {t_2} \cdot 2(t_1^2 - t_2^2) \cdot (2{t_1}{t_2})\}  \cdot [ - sh({t_1}\frac{b}
{2})\sin ({t_2}\frac{b}
{2})]\nonumber\\
&~~~~~~~~~~~~~~~~~~ + \{ {t_1} \cdot 2(t_1^2 - t_2^2) \cdot (2{t_1}{t_2}) + {t_2}[{(t_1^2 - t_2^2)^2} - {(2{t_1}{t_2})^2}]\}  \cdot ch({t_1}\frac{b}
{2})\cos ({t_2}\frac{b}
{2})\} \nonumber\\
&~~~ - \frac{2}
{\varepsilon }\frac{1}
{{{{(t\bar t)}^3}}}\{ \{ (t_1^2 - t_2^2)[{(t_1^2 - t_2^2)^2} - {(2{t_1}{t_2})^2}] - 2(t_1^2 - t_2^2){(2{t_1}{t_2})^2}\}  \cdot [ - ch({t_1}\frac{b}
{2})\sin ({t_2}\frac{b}
{2})]\nonumber\\
&~~~~~~~~~~~~~~~~~ + 2{t_1}{t_2}[3{(t_1^2 - t_2^2)^2} - {(2{t_1}{t_2})^2}]sh({t_1}\frac{b}
{2})\cos ({t_2}\frac{b}
{2})\} \nonumber\\
&~~~ + \frac{2}
{\varepsilon }\frac{1}
{{t\bar t}}\{ (t_1^2 - t_2^2)[ - ch({t_1}\frac{b}
{2})\sin ({t_2}\frac{b}
{2})] + 2{t_1}{t_2}sh({t_1}\frac{b}
{2})\cos ({t_2}\frac{b}
{2})\} \nonumber\\
& ~~~ - [{t_1}\sin ({t_2}b) + {t_2}sh({t_1}b)]\nonumber\\
&~~~ - \frac{2}
{\varepsilon }\frac{b}
{2}\{ {t_1}[ - sh({t_1}\frac{b}
{2})\sin ({t_2}\frac{b}
{2})] + {t_2}ch({t_1}\frac{b}
{2})\cos ({t_2}\frac{b}
{2})\}\nonumber\\
&~~~ - \frac{1}
{\varepsilon }{(\frac{b}
{2})^2}t\bar tch({t_1}\frac{b}
{2})\sin ({t_2}\frac{b}
{2}).
\end{align}
Clearly, $f(0;\varepsilon ) = 0$.

\subsection{Determining ${b_0}$  and ${\varepsilon _0}$  satisfying (\ref{eq2.14}) and (\ref{eq2.15})}

 Let
\begin{equation}\label{eq2.17}
{t_1} = \alpha {t_2},
\end{equation}
then
\begin{equation}\label{eq2.18}
\alpha  > 12
\end{equation}
from ${t_1} > 6$ and $0 < {t_2} < \frac{1}
{2}$, and
\begin{align}\label{eq2.19}
&~~~~~f(b;\varepsilon )\nonumber\\
&=-\frac{1}{\varepsilon}(\frac{b}{2})^2\frac{1}{a^2+1}\{[(\alpha ^2-1)^2 t_2^2-4\alpha ^2t_2^2]\cdot[-ch(\alpha t_2\frac{b}{2})sin (t_2\frac{b}{2})]\nonumber\\
&~~~~~~~~~~~~~~~~~~~~~~~~~ +2(\alpha^2-1)\cdot 2\alpha t_2^2 sh(\alpha t_2 \frac{b}{2})cos(t_2\frac{b}{2})\} \nonumber\\
&~~~+\frac{2}{\varepsilon}\frac{b}{2}\frac{1}{(\alpha^2+1)^2} \{ \{ \alpha t_2[(\alpha^2-1)^2-4\alpha^2]       -t_2\cdot 2(\alpha^2-1)\cdot 2\alpha\}\cdot [-sh(\alpha t_2\frac{b}{2})sin(t_2\frac{b}{2})]            \nonumber\\
&~~~~~~~~~~~~~~~~~~~~~~~~~~+ \left\{\alpha t_2\cdot 2(\alpha^2-1)\cdot 2\alpha +t_2[(\alpha^2-1)^2-4\alpha^2]\right\}\cdot ch(\alpha t_2 \frac{b}{2}) cos(t_2\frac{b}{2})                         \}\nonumber\\
&~~~ - \frac{2}
{\varepsilon }\frac{1}
{{{{({\alpha ^2} + 1)}^3}}}\{ \{ ({\alpha ^2} - 1)[{({\alpha ^2} - 1)^2} - 4{\alpha ^2}] - 2({\alpha ^2} - 1)4{\alpha ^2}\}  \cdot [ - ch(\alpha {t_2}\frac{b}
{2})\sin ({t_2}\frac{b}
{2})]\nonumber\\
&~~~~~~~~~~~~~~~~~~~~~~ + 2\alpha [3{({\alpha ^2} - 1)^2} - 4{\alpha ^2}]sh(\alpha {t_2}\frac{b}
{2})\cos ({t_2}\frac{b}
{2})\}\nonumber\\
&~~~ + \frac{2}
{\varepsilon }\frac{1}
{{{\alpha ^2} + 1}}\{ ({\alpha ^2} - 1)[ - ch(\alpha {t_2}\frac{b}
{2})\sin ({t_2}\frac{b}
{2})] + 2\alpha sh(\alpha {t_2}\frac{b}
{2})\cos ({t_2}\frac{b}
{2})\} \nonumber\\
&~~~ - [\alpha {t_2}\sin ({t_2}b) + {t_2}sh(\alpha {t_2}b)]\nonumber\\
&~~~ - \frac{2}
{\varepsilon }\frac{b}
{2}\{ \alpha {t_2}[ - sh(\alpha {t_2}\frac{b}
{2})\sin ({t_2}\frac{b}
{2})] + {t_2}ch(\alpha {t_2}\frac{b}
{2})\cos ({t_2}\frac{b}
{2})\} \nonumber\\
&~~~ - \frac{1}
{\varepsilon }{(\frac{b}
{2})^2}({\alpha ^2} + 1)t_2^2ch(\alpha {t_2}\frac{b}
{2})\sin ({t_2}\frac{b}
{2})\nonumber\\
&=\frac{1}
{\varepsilon }{(\frac{b}
{2})^2}\frac{{t_2^2}}
{{{\alpha ^2} + 1}}[({\alpha ^4} - 6{\alpha ^2} + 1) ch(\alpha {t_2}\frac{b}
{2})\sin ({t_2}\frac{b}
{2}) - 4\alpha ({\alpha ^2} - 1)sh(\alpha {t_2}\frac{b}
{2})\cos ({t_2}\frac{b}
{2})]\nonumber\\
&~~~ + \frac{2}
{\varepsilon }\frac{b}
{2}\frac{{{t_2}}}
{{{{({\alpha ^2} + 1)}^2}}}\{ ({\alpha ^5} - 10{\alpha ^3} + 5\alpha )[ - sh(\alpha {t_2}\frac{b}
{2})\sin ({t_2}\frac{b}
{2})]\nonumber\\
&~~~~~~~~~~~~~~~~~~~~~~ + (5{\alpha ^4} - 10{\alpha ^2} + 1)ch(\alpha {t_2}\frac{b}
{2})\cos ({t_2}\frac{b}
{2})\} \nonumber\\
&~~~ + \frac{2}
{\varepsilon }\frac{1}
{{{{({\alpha ^2} + 1)}^3}}}\{ [({\alpha ^2} - 1)({\alpha ^4} - 6{\alpha ^2} + 1) - 8{\alpha ^2}({\alpha ^2} - 1)]ch(\alpha {t_2}\frac{b}
{2})\sin ({t_2}\frac{b}
{2})\nonumber\\
&~~~~~~~~~~~~~~~~~~~~~~~ - 2\alpha [3{({\alpha ^2} - 1)^2} - 4{\alpha ^2}]sh(\alpha {t_2}\frac{b}
{2})\cos ({t_2}\frac{b}
{2})\} \nonumber\\
&~~~ + \frac{2}
{\varepsilon }\frac{1}
{{{\alpha ^2} + 1}}\{ ({\alpha ^2} - 1)[ - ch(\alpha {t_2}\frac{b}
{2})\sin ({t_2}\frac{b}
{2})] + 2\alpha sh(\alpha {t_2}\frac{b}
{2})\cos ({t_2}\frac{b}
{2})\} \nonumber\\
&~~~- [\alpha {t_2}\sin ({t_2}b) + {t_2}sh(\alpha {t_2}b)]\nonumber\\
&~~~ + \frac{2}
{\varepsilon }\frac{b}
{2}{t_2}[\alpha sh(\alpha {t_2}\frac{b}
{2})\sin ({t_2}\frac{b}
{2}) - ch(\alpha {t_2}\frac{b}
{2})\cos ({t_2}\frac{b}
{2})]\nonumber\\
&~~~ - \frac{1}
{\varepsilon }{(\frac{b}
{2})^2}({\alpha ^2} + 1)t_2^2ch(\alpha {t_2}\frac{b}
{2})\sin ({t_2}\frac{b}
{2}),
\end{align}

Since in (\ref{eq2.19}) the first and the seventh terms involve the factor $t_2^2$, the second and the sixth terms involve $t_2$, the third, the fourth and fifth terms do not involve $t_2$, we merge the first term and the seventh term, the second term and the sixth term, and the third term and the fourth term, respectively, and obtain
\begin{align}\label{eq2.20}
&~~~~~f(b;\varepsilon )\nonumber\\
&=\frac{1}
{\varepsilon }{(\frac{b}
{2})^2}\frac{{t_2^2}}
{{{\alpha ^2} + 1}}[({\alpha ^4} - 6{\alpha ^2} + 1) ch(\alpha {t_2}\frac{b}
{2})\sin ({t_2}\frac{b}
{2}) - 4\alpha ({\alpha ^2} - 1)sh(\alpha {t_2}\frac{b}
{2})\cos ({t_2}\frac{b}
{2})\nonumber\\
&~~~~~~~~~~~~~~~~~~~~~ - {({\alpha ^2} + 1)^2}ch(\alpha {t_2}\frac{b}
{2})\sin ({t_2}\frac{b}
{2})]\nonumber\\
&~~~ + \frac{2}
{\varepsilon }\frac{b}
{2}\frac{{{t_2}}}
{{{{({\alpha ^2} + 1)}^2}}}[( - {\alpha ^5} + 10{\alpha ^3} - 5\alpha )sh(\alpha {t_2}\frac{b}
{2})\sin ({t_2}\frac{b}
{2}) + (5{\alpha ^4} - 10{\alpha ^2} + 1)ch(\alpha {t_2}\frac{b}
{2})\cos ({t_2}\frac{b}
{2})\nonumber\\
&~~~~~~~~~~~~~~~~~~~~~~~~ + \alpha {({\alpha ^2} + 1)^2}sh(\alpha {t_2}\frac{b}
{2})\sin ({t_2}\frac{b}
{2}) - {({\alpha ^2} + 1)^2}ch(\alpha {t_2}\frac{b}
{2})\cos ({t_2}\frac{b}
{2})]\nonumber\\
&~~~+ \frac{2}
{\varepsilon }\frac{1}
{{{{({\alpha ^2} + 1)}^3}}}[({\alpha ^2} - 1)({\alpha ^4} - 14{\alpha ^2} + 1)ch(\alpha {t_2}\frac{b}
{2})\sin ({t_2}\frac{b}
{2})\nonumber\\
&~~~~~~~~~~~~~~~~~~~~~~~ - 2\alpha [3{({\alpha ^2} - 1)^2} - 4{\alpha ^2}]sh(\alpha {t_2}\frac{b}
{2})\cos ({t_2}\frac{b}
{2})\nonumber\\
&~~~~~~~~~~~~~~~~~~~~~~~ - ({\alpha ^2} - 1){({\alpha ^2} + 1)^2}ch(\alpha {t_2}\frac{b}
{2})\sin ({t_2}\frac{b}
{2}) + 2\alpha {({\alpha ^2} + 1)^2}sh(\alpha {t_2}\frac{b}
{2})\cos ({t_2}\frac{b}
{2})]\nonumber\\
&~~~ - [\alpha {t_2}\sin ({t_2}b) + {t_2}sh(\alpha {t_2}b)]\nonumber\\
&=\frac{2}
{\varepsilon }{(\frac{b}
{2})^2}\frac{{t_2^2}}
{{{\alpha ^2} + 1}}[( - 4{\alpha ^2})ch(\alpha {t_2}\frac{b}
{2})\sin ({t_2}\frac{b}
{2}) + ( - 2{\alpha ^3} + 2\alpha )sh(\alpha {t_2}\frac{b}
{2})\cos ({t_2}\frac{b}
{2})]\nonumber\\
&~~~ + \frac{2}
{\varepsilon }\frac{b}
{2}\frac{{{t_2}}}
{{{{({\alpha ^2} + 1)}^2}}}[(12{\alpha ^3} - 4\alpha )sh(\alpha {t_2}\frac{b}
{2})\sin ({t_2}\frac{b}
{2}) + (4{\alpha ^4} - 12{\alpha ^2})ch(\alpha {t_2}\frac{b}
{2})\cos ({t_2}\frac{b}
{2})]\nonumber\\
&~~~ + \frac{2}
{\varepsilon }\frac{1}
{{{{({\alpha ^2} + 1)}^3}}}[( - 16{\alpha ^4} + 16{\alpha ^2})ch(\alpha {t_2}\frac{b}
{2})\sin ({t_2}\frac{b}
{2}) + ( - 4{\alpha ^5} + 24{\alpha ^3} - 4\alpha )sh(\alpha {t_2}\frac{b}
{2})\cos ({t_2}\frac{b}
{2})]\nonumber\\
&~~~ - [\alpha {t_2}\sin ({t_2}b) + {t_2}sh(\alpha {t_2}b)]\nonumber\\
&:=\frac{2}
{\varepsilon }Q(b) - [\alpha {t_2}\sin ({t_2}b) + {t_2}sh(\alpha {t_2}b)].
\end{align}

Let us turn to prove that there exists ${b_1} > 0$, such that $Q({b_1}) > 0$. Denote
\begin{equation}\label{eq2.21}
\frac{b}
{2}{t_2} = \sigma ,
\end{equation}
where $\sigma  > 0$ will be determined later, then

\begin{align*}
&~~~~~Q(b) = Q(\frac{{2\sigma }}
{{{t_2}}})\\
&=\frac{1}{\alpha^2+1}[(-4\alpha^2\sigma^2)ch(\alpha\sigma)sin\sigma +(-2\alpha^3+2\alpha)\sigma^2sh(\alpha\sigma)cos\sigma]\\
&~~~ + \frac{1}
{{{{({\alpha ^2} + 1)}^2}}}[(12{\alpha ^3} - 4\alpha )\sigma sh(\alpha \sigma )\sin \sigma  + (4{\alpha ^4} - 12{\alpha ^2})\sigma ch(\alpha \sigma )\cos \sigma ]\\
&~~~ + \frac{1}
{{{{({\alpha ^2} + 1)}^3}}}[( - 16{\alpha ^4} + 16{\alpha ^2})ch(\alpha \sigma )\sin \sigma  + ( - 4{\alpha ^5} + 24{\alpha ^3} - 4\alpha )sh(\alpha \sigma )\cos \sigma ]\\
& = \frac{1}
{{{{({\alpha ^2} + 1)}^3}}}[ - 4{\alpha ^2}{({\alpha ^2} + 1)^2}{\sigma ^2}ch(\alpha \sigma )\sin \sigma  + ( - 2{\alpha ^3} + 2\alpha ){({\alpha ^2} + 1)^2}{\sigma ^2}sh(\alpha \sigma )\cos \sigma \\
&~~~~~~~~~~~~~~~~~ + (12{\alpha ^3} - 4\alpha )({\alpha ^2} + 1)\sigma sh(\alpha \sigma )\sin \sigma  + (4{\alpha ^4} - 12{\alpha ^2})({\alpha ^2} + 1)\sigma ch(\alpha \sigma )\cos \sigma \nonumber\\
&~~~~~~~~~~~~~~~~~ + ( - 16{\alpha ^4} + 16{\alpha ^2})ch(\alpha \sigma )\sin \sigma  + ( - 4{\alpha ^5} + 24{\alpha ^3} - 4\alpha )sh(\alpha \sigma )\cos \sigma ]
\end{align*}

\begin{align}\label{eq2.22}
& = \frac{1}
{{{{({\alpha ^2} + 1)}^3}}}\{ \frac{{{e^{\alpha \sigma }}}}
{2}[ - 4{\alpha ^2}{({\alpha ^2} + 1)^2}{\sigma ^2} + (12{\alpha ^3} - 4\alpha )({\alpha ^2} + 1)\sigma  + ( - 16{\alpha ^4} + 16{\alpha ^2})]\sin \sigma \nonumber\\
&~~~~~~~~~~~~~~~~~~ + \frac{{{e^{\alpha \sigma }}}}
{2}[( - 2{\alpha ^3} + 2\alpha ){({\alpha ^2} + 1)^2}{\sigma ^2}{\text{ + }}(4{\alpha ^4} - 12{\alpha ^2})({\alpha ^2} + 1)\sigma {\text{ + }}( - 4{\alpha ^5} + 24{\alpha ^3} - 4\alpha )]\cos \sigma \nonumber\\
&~~~~~~~~~~~~~~~~~~ +\frac{{{e^{ - \alpha \sigma }}}}
{2}[ - 4{\alpha ^2}{({\alpha ^2} + 1)^2}{\sigma ^2} + ( - 12{\alpha ^3} + 4\alpha )({\alpha ^2} + 1)\sigma  + ( - 16{\alpha ^4} + 16{\alpha ^2})]\sin \sigma \nonumber\\
 &~~~~~~~~~~~~~~~~~~+ \frac{{{e^{ - \alpha \sigma }}}}
{2}[(2{\alpha ^3} - 2\alpha ){({\alpha ^2} + 1)^2}{\sigma ^2} + (4{\alpha ^4} - 12{\alpha ^2})({\alpha ^2} + 1)\sigma  + (4{\alpha ^5} - 24{\alpha ^3} + 4\alpha )]\cos \sigma \} \nonumber\\
&: = \frac{1}
{{{{({\alpha ^2} + 1)}^3}}}\{ \frac{{{e^{\alpha \sigma }}}}
{2}[{g_1}(\sigma )\sin \sigma  + {g_2}(\sigma )\cos \sigma ] + \frac{{{e^{ - \alpha \sigma }}}}
{2}[{g_3}(\sigma )\sin \sigma  + {g_4}(\sigma )\cos \sigma ]\} .
\end{align}
To ${g_1}(\sigma )$  ($\sigma >0$), we see that its image is open side down because its leading coefficient is negative; its discriminant satisfies by using $\alpha >12$  that
\[{(12{\alpha ^3} - 4\alpha )^2}{({\alpha ^2} + 1)^2} - 4[ - 4{\alpha ^2}{({\alpha ^2} + 1)^2}]( - 16{\alpha ^4} + 16{\alpha ^2}) < 0,\]
which gives that for any $\sigma  > 0$,
$$g_1(\sigma)<0.$$
Similarly, it knows that for any $\sigma  > 0$,
$${g_2}(\sigma ) < 0, {g_3}(\sigma ) < 0, {g_4}(\sigma ) > 0.$$
Taking $\sigma  = \frac{{3\pi }}
{2},$ it yields
\[{b_1} = \frac{{3\pi }}
{{{t_2}}},\]
 and by applying $\sin \frac{{3\pi }}
{2} =  - 1$ and $\cos \frac{{3\pi }}
{2} = 0$  that
\[{g_1}(\frac{{3\pi }}
{2})\sin \frac{{3\pi }}
{2} > 0, {g_2}(\frac{{3\pi }}
{2})\cos \frac{{3\pi }}
{2} = 0, {g_3}(\frac{{3\pi }}
{2})\sin \frac{{3\pi }}
{2} > 0, {g_4}(\frac{{3\pi }}
{2})\cos \frac{{3\pi }}
{2} = 0.\]
We reach $Q({b_1}) > 0$  from (\ref{eq2.22}).

On the other hand, take
\[{b_2} = \frac{{5\pi }}
{{{t_2}}},\]
then  $\sin \frac{{5\pi }}
{2} = 1, \cos \frac{{5\pi }}
{2} = 0,$  and
\[{g_1}(\frac{{5\pi }}
{2})\sin \frac{{5\pi }}
{2} < 0,~{g_2}(\frac{{5\pi }}
{2})\cos \frac{{5\pi }}
{2} = 0,~{g_3}(\frac{{5\pi }}
{2})\sin \frac{{5\pi }}
{2} < 0,~{g_4}(\frac{{5\pi }}
{2})\cos \frac{{5\pi }}
{2} = 0.\]
It gets $Q({b_2}) < 0$  from (\ref{eq2.20}).

Now we pick a small ${\varepsilon _1} > 0$ in (\ref{eq2.20}) to arrive at
\begin{equation}\label{eq2.23}
f({b_1};{\varepsilon _1}) > 0, f({b_2};{\varepsilon _1}) < 0.
\end{equation}
\begin{remark}\label{Re2.1}
 One can not obtain directly ${b_0} \in ({b_1},{b_2})$  such that $f({b_0};{\varepsilon _1}) = 0$  for the continuous function $f(b;{\varepsilon _1})$  by using (\ref{eq2.23}) and the intermediate value theorem for continuous functions, because  $b_0$ and $\varepsilon_1$  are difficult to ensure (\ref{eq2.15}). It means that a more careful analysis to derive (\ref{eq2.14}) and (\ref{eq2.15}) is needed.
 \end{remark}

Denote
\begin{equation}\label{eq2.24}
h(b,\varepsilon ) = \int_0^b {v\bar v} dy - \frac{{t\bar t}}
{{10{\varepsilon ^2}}}{(\frac{b}
{2})^2}, \varepsilon >0.
\end{equation}
A direct calculation shows
\begin{equation}\label{eq2.25}
h(b,\varepsilon ) = F(b) + \frac{1}
{\varepsilon }G(b), b \in [{b_1},{b_2}],
\end{equation}
where
\begin{equation}\label{eq2.26}
F(b){\text{ = }}\frac{1}
{4}\left[ {\frac{{{e^{\alpha {t_2}b}} - {e^{ - \alpha {t_2}b}}}}
{{\alpha {t_2}}} + \frac{1}
{{{t_2}}}\sin ({t_2}b)} \right] > 0,
\end{equation}
\begin{align}\label{eq2.27}
G(b)& = {t_1}\int_0^b {\cos [{t_2}(y - \frac{b}
{2})]} ch[{t_1}(y - \frac{b}
{2})]{(y - \frac{b}
{2})^2}dy\nonumber\\
&~~~+ {t_2}\int_0^b {\sin [{t_2}(y - \frac{b}
{2})]} sh[{t_1}(y - \frac{b}
{2})]{(y - \frac{b}
{2})^2}dy.
\end{align}

We now apply (\ref{eq2.25})-(\ref{eq2.27}) to prove (\ref{eq2.14})-(\ref{eq2.15}) and the proofs of (\ref{eq2.25})-(\ref{eq2.27}) will be given in Section 3.

Note that $f(b,\varepsilon ),~Q(b)$  and $G(b)$  are all continuous on $[{b_1}, {b_2}]$ and may change sign; their zeros are isolated. Let us prove that there exist ${b_0} \in ({b_1},{b_2})$  and ${\varepsilon _0} > 0$, such that
\begin{equation}\label{eq2.28}
f({b_0},{\varepsilon _0}) = 0, F({b_0}) + \frac{1}
{{{\varepsilon _0}}}G({b_0}) \ne 0.
\end{equation}
We will treat the case that $G(b)$  has at most a zero in $[{b_1},{b_2}]$; other case can be inspected similarly.

(1)  If $G(b)$  has no any zero on $[{b_1},{b_2}]$, then $G(b)$  is always positive or
negative on $[{b_1},{b_2}]$, and there is  $c_0>0$ such that
\begin{equation}\label{eq2.29}
G(b) \geqslant {c_0} \text{~or~} G(b) \leqslant  - {c_0}, b \in [{b_1},{b_2}].
\end{equation}

If $G(b) \geqslant {c_0}, b \in [{b_1},{b_2}],$
 then by noting $F(b)>0$  in (\ref{eq2.26}) that for any $\varepsilon >0$  and $b \in ({b_1},{b_2})$, it holds $F(b) + \frac{1}
{\varepsilon }G(b) > 0,$  and hence for a specific ${\varepsilon _1} > 0$,
\[F(b) + \frac{1}
{{{\varepsilon _1}}}G(b) \ne 0.\]
Using (\ref{eq2.23}), the intermediate value theorem for continuous functions and the fact that $f(b,{\varepsilon _1})$ is continuous on $[{b_1},{b_2}]$, there exists ${b_0} \in ({b_1},{b_2})$, such that $f({b_0},{\varepsilon _1}) = 0$. The previous formula also implies $F({b_0}) + \frac{1}
{{{\varepsilon _1}}}G({b_0}) \ne 0$ and (\ref{eq2.28}) is proved.

If $G(b) \leqslant  - {c_0}$, $b \in [{b_1},{b_2}]$, then we take ${\varepsilon _2} = \frac{{{c_0}}}
{{\mathop {\max }\limits_{b \in [{b_1},{b_2}]} F(b)}} > 0$  by using (\ref{eq2.26}) and have
\[F(b) + \frac{1}
{{{\varepsilon _2}}}G(b) \leqslant \mathop {\max }\limits_{b \in [{b_1},{b_2}]} F(b) + \frac{{ - {c_0}}}
{{{\varepsilon _2}}} = 0\]
and
\begin{equation}\label{eq2.30}
F(b) + \frac{1}
{\varepsilon }G(b) < F(b) + \frac{1}
{{{\varepsilon _2}}}G(b) \leqslant 0, \varepsilon <\varepsilon _2.
\end{equation}
By (\ref{eq2.23}) and (\ref{eq2.20}), it gives
\[f({b_1},\varepsilon ) > 0,f({b_2},\varepsilon ) < 0, \text{~for any~}\varepsilon  < {\varepsilon _1}.\]
Choosing  ${\varepsilon _0}, 0 < {\varepsilon _0} < \min ({\varepsilon _1},{\varepsilon _2}),$ for the function $f(b,{\varepsilon _0})$, we see from the above conclusions that there exists ${b_0} \in ({b_1},{b_2})$, such that $f({b_0},{\varepsilon _0}) = 0$; and $F({b_0}) + \frac{1}
{{{\varepsilon _0}}}G({b_0}) \ne 0$ by (\ref{eq2.30}). Now (\ref{eq2.28}) is verified.

(2)  If  $G(b)$  has a zero $b'$ on  $[{b_1},{b_2}]$, then it happens three possibilities:
(i) $b' \in ({b_1},{b_2})$; (ii) $b'{\text{ = }}{b_1}$; (iii) $b'{\text{ = }}{b_2}$. Let us deal with them in the sequel.

(i) If $G(b') = 0$, $b' \in ({b_1},{b_2})$, then we consider two cases for the continuous function $Q(b): Q(b') > 0$ and $Q(b') \leqslant 0$.

To the case $Q(b') > 0$, we know from keeping sign property of continuous functions that there is ${b'_1} > b'$ (so ${b_1} < b' < {b'_1} < {b_2}$), such that $Q({b'_1}) > 0$; inexpensively, one can select ${\varepsilon _1} > 0$, such that $f({b'_1},{\varepsilon _1}) > 0,f({b_2},{\varepsilon _1}) < 0$. Since $G(b)$ is always positive or negative on $[{b'_1},{b_2}]$, we can use the same way in the case (1) to derive that there exist    ${\varepsilon _0}, 0 < {\varepsilon _0} < {\varepsilon _1}$ and  ${b_0} \in [{b'_1},{b_2}] \subset [{b_1},{b_2}]$, such that $f({b_0},{\varepsilon _0}) = 0$ and $F({b_0}) + \frac{1}
{{{\varepsilon _0}}}G({b_0}) \ne 0.$ Hence (\ref{eq2.28}) is proved.

To the case $Q(b') \le 0$, we consider $Q(b') < 0$ and $Q(b') = 0$, respectively. If $Q(b') < 0$, then by the fact that $Q(b)$ is continuous, there is ${b'_2} < b',$ such that $Q({b'_2}) < 0,$ so $f({b'_2},{\varepsilon _1}) = \frac{2}{{{\varepsilon _1}}}Q({b'_2}) - I({b'_2}) < 0,$ where $ I({b'_2})=\alpha {t_2} sin({t_2}{b'_2})+{t_2} sh(\alpha{t_2}{b'_2})$, ${\varepsilon _1}$ sees \eqref{eq2.23}. It implies that for $0 < \varepsilon < {\varepsilon _1}$,
\[f({b'_2},\varepsilon ) = \frac{2}{\varepsilon }Q({b'_2}) - I({b'_2}) < \frac{2}{{{\varepsilon _1}}}Q({b'_2}) - I({b'_2}) < 0,\]
and it holds from $Q({b_1}) > 0,$
\[f({b_1},\varepsilon ) = \frac{2}{\varepsilon }Q({b_1}) - I({b_1}) > \frac{2}{{{\varepsilon _1}}}Q({b_1}) - I({b_1}) > 0.\]
Now $G(b)$ is always positive or negative on $[{b_1},{b'_2}]$, and we can prove \eqref{eq2.28} by the same way in the case (1).

If $Q(b') = 0$, then for any $0 < \varepsilon  < {\varepsilon _1},$ it follows
\[f(b',\varepsilon ) = \frac{2}{\varepsilon }Q(b') - I(b') =  - I(b') < 0;\]
combining $f({b_1},\varepsilon ) > 0,$ there exists ${b_0} \in ({b_1},b'),$ such that
\[f({b_0},\varepsilon ) = 0\;(0 < \varepsilon  < {\varepsilon _1}).\]
Noting $G({b_0}) \ne 0,$ we see that when $G({b_0}) > 0$, it holds
\[F({b_0}) + \frac{1}{\varepsilon }G({b_0}) > 0\; (0 < \varepsilon  < {\varepsilon _1}),\]
so \eqref{eq2.28} is proved; when $G({b_0}) < 0$, we can take ${\varepsilon _0}$ satisfying
\[0 < {\varepsilon _0} < \min \left\{ { - \frac{{G({b_0})}}{{F({b_0})}},{\varepsilon _1}} \right\}\]
to obtain $F({b_0}) + \frac{1}{{{\varepsilon _0}}}G({b_0}) < 0$ and $f({b_0},{\varepsilon _0}) = 0$ from $f({b_0},\varepsilon ) = 0$
($0 < \varepsilon  < {\varepsilon _1}$), thus \eqref{eq2.28} is proved.

(ii) If $G({b_1}) = 0$, then from  $f({b_1},{\varepsilon _1}) > 0$ (see (\ref{eq2.23})) and keeping sign property of the continuous function $f(b,{\varepsilon _1})$, there exists ${b'_1},$  ${b_1} < {b'_1} < {b_2},$  such that $f({b'_1},{\varepsilon _1}) > 0,$  $f({b_2},{\varepsilon _1}) < 0$, and $G(b)$ is always positive or negative on $[{b'_1},{b_2}]$. Now we argue as in (1).

(iii) If $G({b_2}) = 0$, then from $f({b_2},{\varepsilon _1}) < 0$ (see (\ref{eq2.23})), there is ${b'_2},$ ${b_1} < {b'_2} < {b_2},$ such that $f({b'_2},{\varepsilon _1}) < 0,f({b_1},{\varepsilon _1}) > 0,$  and $G(b)$  is always positive or negative on $[{b_1},{b'_2}]$. We return again to (1).

At this point, (\ref{eq2.28}) is derived. Since (\ref{eq2.28}) indicates (\ref{eq2.14}) and (\ref{eq2.15}), we obtain from (\ref{eq2.13}) that  $t^2$ is real.

\subsection{ Finishing the proof of Theorem 1 }

From
\[{t^2} = ({t_1} + i{t_2}) \cdot ({t_1} + i{t_2}) = t_1^2 - t_2^2 + i \cdot 2{t_1}{t_2}\]
and the claim that $t^2$  is real, it infers ${t_1}{t_2} = 0$. Noting  ${t_2} \ne 0$ and $t_1>6$, it is impossible, therefore, ${t_2} = 0$, i.e., (\ref{eq1.6}) is true. The proof of Theorem 1 is ended.

\section{Proofs of (\ref{eq2.25})-(\ref{eq2.27}) }
In $h(b,\varepsilon ) = \int_0^b {v\bar v} dy - \frac{{t\bar t}}
{{10{\varepsilon ^2}}}{(\frac{b}
{2})^2},\varepsilon  > 0,$  we note via (\ref{eq2.3}) that
\[v = \frac{{{e^{t(y - \frac{b}
{2})}} + {e^{ - t(y - \frac{b}
{2})}}}}
{2} + \frac{1}
{{2\varepsilon }}t{(y - \frac{b}
{2})^2},\]
and use
\[t = {t_1} + i{t_2},\]
\[{e^{t(y - \frac{b}
{2})}} = {e^{{t_1}(y - \frac{b}
{2})}}\left[ {\cos ({t_2}(y - \frac{b}
{2})) + i\sin ({t_2}(y - \frac{b}
{2}))} \right],\]
\[{e^{ - t(y - \frac{b}
{2})}} = {e^{ - {t_1}(y - \frac{b}
{2})}}\left[ {\cos ({t_2}(y - \frac{b}
{2})) - i\sin ({t_2}(y - \frac{b}
{2}))} \right],\]
\[{e^{t(y - \frac{b}
{2})}} + {e^{ - t(y - \frac{b}
{2})}} = 2\cos ({t_2}(y - \frac{b}
{2}))ch({t_1}(y - \frac{b}
{2})) + 2i\sin ({t_2}(y - \frac{b}
{2}))sh({t_1}(y - \frac{b}
{2}))\]
to know
\[v = \cos ({t_2}(y - \frac{b}
{2}))ch({t_1}(y - \frac{b}
{2})) + \frac{{{t_1}}}
{{2\varepsilon }}{(y - \frac{b}
{2})^2} + i\left[ {sin({t_2}(y - \frac{b}
{2}))sh({t_1}(y - \frac{b}
{2})) + \frac{{{t_2}}}
{{2\varepsilon }}{(y - \frac{b}
{2})^2}} \right],\]
\[\bar v = \cos ({t_2}(y - \frac{b}
{2}))ch({t_1}(y - \frac{b}
{2})) + \frac{{{t_1}}}
{{2\varepsilon }}{(y - \frac{b}
{2})^2} - i\left[ {sin({t_2}(y - \frac{b}
{2}))sh({t_1}(y - \frac{b}
{2})) + \frac{{{t_2}}}
{{2\varepsilon }}{(y - \frac{b}
{2})^2}} \right],\]
and
\begin{align*}
v\bar v &= {\cos ^2}({t_2}(y - \frac{b}
{2}))c{h^2}({t_1}(y - \frac{b}
{2})) + {\sin ^2}({t_2}(y - \frac{b}
{2}))s{h^2}({t_1}(y - \frac{b}
{2}))\\
&~~~+ \frac{{{t_1}}}
{\varepsilon }\cos ({t_2}(y - \frac{b}
{2}))ch({t_1}(y - \frac{b}
{2})){(y - \frac{b}
{2})^2} + \frac{{{t_2}}}
{\varepsilon }\sin ({t_2}(y - \frac{b}
{2}))sh({t_1}(y - \frac{b}
{2})){(y - \frac{b}
{2})^2}\\
&~~~ + \frac{{t_1^2}}
{{4{\varepsilon ^2}}}{(y - \frac{b}
{2})^4} + \frac{{t_2^2}}
{{4{\varepsilon ^2}}}{(y - \frac{b}
{2})^4}.
\end{align*}
Since
\[\int_0^b {\frac{{t_1^2 + t_2^2}}
{{4{\varepsilon ^2}}}} {(y - \frac{b}
{2})^4}dy = \frac{{t\bar t}}
{{10{\varepsilon ^2}}}{(\frac{b}
{2})^5},\]
it follows
\begin{align*}
&~~~~~h(b,\varepsilon ) \\
&= \int_0^b {{{\cos }^2}({t_2}(y - \frac{b}
{2}))} c{h^2}({t_1}(y - \frac{b}
{2}))dy + \int_0^b {{{\sin }^2}({t_2}(y - \frac{b}
{2}))s{h^2}({t_1}(y - \frac{b}
{2}))} dy\\
&~~~ + \frac{1}
{\varepsilon }\left[ {{t_1}\int_0^b {\cos ({t_2}(y - \frac{b}
{2}))} ch({t_1}(y - \frac{b}
{2})){(y - \frac{b}
{2})^2}dy + {t_2}\int_0^b {\sin ({t_2}(y - \frac{b}
{2}))} sh({t_1}(y - \frac{b}
{2})){(y - \frac{b}
{2})^2}dy} \right],
\end{align*}
which is (\ref{eq2.25}), the first term of the right hand side is $F(b)$  and the parenthesis in the second term is $G(b)$, i.e., (\ref{eq2.27}). Noting
\begin{align*}
&~~~~~\int_0^b {{{\cos }^2}({t_2}(y - \frac{b}
{2}))} c{h^2}({t_1}(y - \frac{b}
{2}))dy\\
& = \int_0^b {\frac{{1 + \cos (2{t_2}(y - \frac{b}
{2}))}}
{2}}  \cdot \frac{{{e^{2{t_2}(y - \frac{b}
{2})}} + 2 + {e^{ - 2{t_2}(y - \frac{b}
{2})}}}}
{2}dy\\
& = \frac{1}
{8}\int_0^b {\left[ {{e^{2{t_2}(y - \frac{b}
{2})}} + 2 + {e^{ - 2{t_2}(y - \frac{b}
{2})}} + \cos (2{t_2}(y - \frac{b}
{2})){e^{2{t_2}(y - \frac{b}
{2})}}} \right.} \\
&~~~\left. { + 2\cos (2{t_2}(y - \frac{b}
{2})) + \cos (2{t_2}(y - \frac{b}
{2})){e^{ - 2{t_2}(y - \frac{b}
{2})}}} \right]dy,
\end{align*}
\begin{align*}
&~~~\int_0^b {{{\sin }^2}({t_2}(y - \frac{b}
{2}))} s{h^2}({t_1}(y - \frac{b}
{2}))dy\\
& = \int_0^b {\frac{{1 - \cos (2{t_2}(y - \frac{b}
{2}))}}
{2}}  \cdot \frac{{{e^{2{t_2}(y - \frac{b}
{2})}} - 2 + {e^{ - 2{t_2}(y - \frac{b}
{2})}}}}
{2}dy\\
& = \frac{1}
{8}\int_0^b {\left[ {{e^{2{t_2}(y - \frac{b}
{2})}} - 2 + {e^{ - 2{t_2}(y - \frac{b}
{2})}} - \cos (2{t_2}(y - \frac{b}
{2})){e^{2{t_2}(y - \frac{b}
{2})}}} \right.} \\
&~~~\left. { + 2\cos (2{t_2}(y - \frac{b}
{2})) - \cos (2{t_2}(y - \frac{b}
{2})){e^{ - 2{t_2}(y - \frac{b}
{2})}}} \right]dy,
\end{align*}
we have
\begin{align*}
F(b)& = \int_0^b {{{\cos }^2}({t_2}(y - \frac{b}
{2}))} c{h^2}({t_1}(y - \frac{b}
{2}))dy + \int_0^b {{{\sin }^2}({t_2}(y - \frac{b}
{2}))} s{h^2}({t_1}(y - \frac{b}
{2}))dy\\
& = \frac{1}
{8}\int_0^b {\left[ {2{e^{2{t_1}(y - \frac{b}
{2})}} + 2{e^{ - 2{t_1}(y - \frac{b}
{2})}} + 4\cos (2{t_2}(y - \frac{b}
{2}))} \right]} dy\\
& = \frac{1}
{8}\left[ {\frac{2}
{{{t_1}}}{e^{2{t_1}\frac{b}
{2}}} - \frac{2}
{{{t_1}}}{e^{ - 2{t_1}\frac{b}
{2}}} + \frac{4}
{{{t_2}}}\sin (2{t_2}\frac{b}
{2})} \right]\\
& = \frac{1}
{4}\left[ {\frac{{{e^{{t_1}b}} - {e^{ - {t_1}b}}}}
{{{t_1}}} + \frac{2}
{{{t_2}}}\sin ({t_2}b)} \right]\\
& = \frac{1}
{4}\left[ {\frac{{{e^{\alpha {t_2}b}} - {e^{ - \alpha {t_2}b}}}}
{{\alpha {t_2}}} + \frac{2}
{{{t_2}}}\sin ({t_2}b)} \right],
\end{align*}
where ${t_1} = \alpha {t_2}$  is applied in the last equality.

Let us prove $F(b) > 0, b \in [{b_1},{b_2}].$    By ${e^{\alpha {t_2}b}} > 1 + \alpha {t_2}b, {e^{ - \alpha {t_2}b}} < 1, \sin ({t_2}b) \geqslant  - 1, b{t_2} \in [3\pi ,5\pi ],$ it implies
\[F(b) = \frac{1}
{4}\left[ {\frac{{{e^{\alpha {t_2}b}} - {e^{ - \alpha {t_2}b}}}}
{{\alpha {t_2}}} + \frac{2}
{{{t_2}}}\sin ({t_2}b)} \right] > \frac{1}
{4}\left( {\frac{{1 + \alpha {t_2}b - 1}}
{{\alpha {t_2}}} - \frac{2}
{{{t_2}}}} \right) = \frac{1}
{4}\left( {b - \frac{2}
{{{t_2}}}} \right) \geqslant \frac{1}
{4}\left( {\frac{{3\pi }}
{{{t_2}}} - \frac{2}
{{{t_2}}}} \right) > 0.\]
Now (\ref{eq2.26}) is proved.

\textbf{Acknowledgement}

I am grateful to Qianqiao Guo, Xueli Bai, J. G\'{e}linas, Wenjuan Li, Huiju Wang, Xiaoxue Ji and Shihong Zhang for their comments, to Junqiang Han and Pengfei Ma for their help in numeric analysis and to Leyun Wu for her type setting.

 \textbf{Conflict of Interest Statement}

 There is no conflict of interest.

\textbf{Data availability statement}

No data was used.

%%appendix------------------------------------------------------------------------

%%-------------------------the followings are the proofs of e102, e103

%% bibliography--------------------------------------------------------------------

\end{document}